\documentclass[preprint,12pt]{article}
\usepackage{amssymb}
\usepackage{mathrsfs}
\usepackage{amsfonts}
\usepackage{graphicx}
\usepackage{colortbl,dcolumn}
\usepackage{amsmath}
\usepackage{psfrag}
\usepackage{booktabs}
\usepackage{cite}
\usepackage{verbatim}

\usepackage{amsthm}
\numberwithin{equation}{section}
\usepackage[top=0.8in, bottom=0.8in, left=0.8in, right=0.8in]{geometry}

\newcommand{\R}{\mathbb{R}}
\newcommand{\N}{\mathbb{N}}
\newcommand{\E}{\mathbb{E}}
\renewcommand{\P}{\mathbb{P}}

\newcommand{\dd}{\text{d}}

\newtheorem{thm}{Theorem}[section]

\newtheorem{lem}[thm]{Lemma}
\newtheorem{prop}[thm]{Proposition}
\newtheorem{cor}[thm]{Corollary}

\begin{document}
\title{Multilevel Monte Carlo methods for 
positivity-preserving approximations of 
the Heston 3/2-model
%Improved split-step implicit Milstein type methods for SDEs with non-globally Lipschitz diffusion coefficients
%New implicit Milstein type methods for SDEs with non-Lipschitz coefficients: applications to financial models
%Implicit Milstein type methods for SDEs with non-globally Lipschitz coefficients
%and their applications of approximating financial models
%Mean-square convergence rates of implicit Milstein methods for SDEs with super-linearly growing coefficients
%and their applications to stochastic volatility processes
 \footnotemark[1]}

\author{
Xiaojuan Wu$\,^{\text{a,b}}$, Siqing Gan$^{\,\text{b}}$\footnotemark[2]
%,  
\\
\footnotesize $\,^\text{a}$ School of Mathematics and Statistics, Hunan First Normal University, Changsha, China
\\
\footnotesize $\,^\text{b}$ School of Mathematics and Statistics, HNP-LAMA, Central South University, Changsha, China
}

      \date{}
       \maketitle

       \footnotetext{\footnotemark[1] This work was supported by Natural Science Foundation of China (12371417, 11971488).
                }
        \footnotetext{\footnotemark[2] Corresponding author: sqgan@csu.edu.cn}
%        \footnotetext{\footnotemark[3] ...@...}
%        \footnotetext{\footnotemark[4] ...@...}
       \begin{abstract}
          {\rm\small  
          This article is concerned with the
          multilevel Monte Carlo (MLMC) methods for approximating expectations of some functions of the solution to the Heston 3/2-model 
          from mathematical finance, 
          which takes values in $(0, \infty)$
          and possesses superlinearly growing drift and diffusion coefficients. 
          To discretize the SDE model, a new Milstein-type scheme is proposed to produce independent sample paths.
          The proposed scheme can be explicitly solved and is positivity-preserving unconditionally, i.e., for any time step-size $h >0$. This  positivity-preserving property for large discretization time steps is particularly desirable in the MLMC setting.
          Furthermore, a mean-square convergence rate of
          order one is proved in the non-globally Lipschitz regime, which is not trivial, as the diffusion coefficient grows super-linearly. 
          The obtained order-one convergence in turn promises the desired relevant variance of the multilevel estimator
        and justifies the optimal complexity $\mathcal{O}(\epsilon^{-2})$ for the MLMC approach, where
        $\epsilon > 0$ is the required target accuracy. 
        Numerical experiments are finally reported to 
        confirm the theoretical findings.
          %covering SDEs  with super-linearly growing diffusion coefficients, such as the popular 3/2-model .
           } \\

%\textbf{PACS: 02.60.Lj}

\end{abstract}
%%%%%%%
%%%%%%

%%%%%
%%%%%
%
\section{Introduction}	
The pricing of financial assets plays an important role 
in mathematical finance \cite{glasserman2004monte,higham2004introduction}, where stochastic differential equations (SDEs) 
\begin{equation} \label{intro-eq:SODE}
%\left\{
%    \begin{array}{ll}
    \dd X(t)  = f ( X (t) ) \, \dd t + g ( X (t) ) \,\dd W (t),
    \quad
    t \in (0, T],
    \quad
    X_0 = x_0
% \end{array}\right.
\end{equation}
are usually used to capture the behaviors of financial systems.
More accurately, there has been 
an increasing interest in quantification of expectations of some functions of the solution $X({T})$ to SDEs, that is, one wants to compute
\begin{equation} \label{introduction:expectation-of-functional-of-sde}
    \mathbb{E} 
    \left[
    \phi\left( X ({T}) \right) 
    \right]
\end{equation}
for some given (payoff) function $\phi$. Usually,
the Monte Carlo (MC) approach, in combination with the Euler-Maruyama (EM) time-stepping scheme, is adopted for this quantification, but at a high computational expense.
Specifically, the computational cost of approximating the quantity \eqref{introduction:expectation-of-functional-of-sde} would be 
$\mathcal{O}(\epsilon^{-3})$, 
where $\epsilon > 0$ is the required target accuracy
measured by the root-mean-square error \cite{duffie1995efficient}.

Recently, Giles developed a
multilevel Monte Carlo (MLMC) method \cite{giles2008multilevel} to significantly reduce
the computational cost of approximating the quantity \eqref{introduction:expectation-of-functional-of-sde}.
Different from the single level MC method, 
the MLMC method relies on a sequence of levels of MC approximations, produced by a strongly convergent time-stepping scheme used
to discretize SDEs \eqref{intro-eq:SODE} with different timestep sizes (see Section \ref{Milstein-sec:conv-rates-global-polynomial} for more details).
The key idea of the MLMC to reduce the overall computational complexity lies in the fact that more paths are simulated on the coarser level and fewer paths are simulated on the finer level \cite{giles2008multilevel}. This way promises a significant reduction in the overall computational cost to $\mathcal{O}(\epsilon ^{-2} (\log(1/\epsilon))^{2})$ when combined with the EM discretization scheme, and to $\mathcal{O}(\epsilon ^{-2})$ when 
combined with the Milstein time-stepping scheme \cite{giles2008improved}. 

It is worthwhile to mention that, the above theoretical
results were obtained for SDEs with globally Lipschitz
continuous coefficients. What if the globally Lipschitz condition is violated? As answered by \cite{hutzenthaler2011strong}, the popularly used EM method fails to converge in both weak and strong senses, when used to solve a large class of SDEs with super-linearly growing coefficients. Furthermore,  Hutzenthaler and Jentzen \cite{hutzenthaler2013divergence} showed that,
the MLMC approach combined with the EM method 
also fails to converge to the desired quantity 
\eqref{introduction:expectation-of-functional-of-sde}, 
when the considered SDEs have super-linearly growing coefficients.  To overcome the  issue, the authors of \cite{hutzenthaler2013divergence} introduced 
a tamed Euler time discretization scheme and obtained the convergence in the multilevel setting. 
Also, the authors of \cite{guo2018multi}
embedded the MLMC method with a truncated EM method \cite{mao2015truncated} to
handle SDEs with non-globally coefficients.

In this article, we focus on the quantification of
\eqref{introduction:expectation-of-functional-of-sde} for
the Heston $3/2$-model \cite{heston1997simple,ahn1999parametric}:
%\cite{lewis2000option}:
%
\begin{equation}
\label{eq:intro-32model}
\begin{split}
    d X (t) = X(t) ( \mu - \alpha X(t) ) dt + \beta X(t) ^ {3/2} \, \dd W(t), \quad
 \mu , \alpha , \beta > 0,
 \quad
     X(0) = X_0 > 0.
\end{split}
\end{equation}
Such an equation can be also used for modelling term structure dynamics \cite{ahn1999parametric},
and can be regarded as a stochastic extension to the logistic equation \cite{gard1988introduction}.
To approximate \eqref{eq:intro-32model}, 
the commonly used EM scheme is 
not a good candidate for two reasons. On the one hand,
the EM is not able to produce convergent approximations,
as both drift and diffusion coefficients of the model grow polynomially. On the other hand, the EM scheme
fails to preserve the positivity of the 
original financial model, which is particularly desirable for the scheme to be well-defined.
To address the aforementioned issues, 
Higham, Mao and Szpruch \cite{higham2013convergence} constructed a class of
positivity preserving schemes for the model
\eqref{eq:intro-32model}, but without a convergence rate revealed.
In this paper we propose
a new Milstein-type scheme \eqref{Chpter5Milstein-eq:scheme} (see also \eqref{Chpter5Milstein-eq:scheme2} and \eqref{eq:Spli-step-doub-implicit-Mil}), which can be regarded as an improved version of the split-step Milstein method introduced in \cite{wu2022split-step}.
The proposed scheme applied to the underlying model 
can be explicitly solved, as shown by
\eqref{eq:Spli-step-doub-implicit-Mil}, and
is proved to be positivity-preserving 
for any time step-size $h >0$.
This unconditional positivity-preserving property  is particularly desirable as large discretization time steps are naturally used in the MLMC setting. 

It is worthwhile to point out that the  
Heston $3/2$-model has
super-linearly growing coefficients, which raises many new difficulties in the error analysis.
To overcome such difficulties, we borrow the idea in \cite{wu2022split-step} 
and introduce an auxiliary process
$\widetilde{X} ( t_n )$ defined by \eqref{eq:Chapter5Aux-X}.
In this way, a mean-square convergence rate of
order one is proved in the non-globally Lipschitz regime
(cf. Theorem \ref{thm:Chapter5Milstein-convergence-rate-32model}),
without relying on a priori estimates for the high-order moment of numerical approximations. This in turn promises the desired relevant variance of the multilevel estimator (see Lemma \ref{lem:variance})
        and justifies the optimal complexity $\mathcal{O}(\epsilon^{-2})$ for the MLMC approach (see Corollary \ref{cor:optimal-complexity}). 

%\cite{wang2020mean,higham2015introduction}

The rest of this article is organized as follows. 
The next section gives the proposed scheme and formulates some assumptions. 
The optimal complexity of the MLMC approach
combined with the proposed time-stepping scheme
is obtained in section \ref{Milstein-sec:conv-rates-global-polynomial}.
Finally, some numerical simulations are provided  in section \ref{Milstein-sect:numer-results} and a short summary is made in section \ref{sec:conclusion}.

%In the past decade, an intensive study of ... 
%In the majority of these results it is assumed that...

%stochastic theta methods \cite{wang2020mean}
%\cite{sabanis2016euler}

%This work is concerned with mean-square approximations of SDEs with non-Lipschitz coefficients.
           
%          Strong convergence rates of SST methods with 
%          parameters $\theta \in [\tfrac12, 1]$ were derived in \cite{yue2014strong} in a non-Lipschitz setting,  
%          while the diffusion coefficients are required to be globally Lipschitz continuous. 

\section{The proposed time discretization scheme}
%\section{Preliminaries}
\label{Milstein-sect:setting}
%

%Let $ d,m, N \in \mathbb{N} $, $ T \in (0, \infty) $ and let $ \left( \Omega, \mathcal{ F }, \P,  \{ \mathcal{ F }_t \}_{ t \in [0,T] }
%\right) $ be a filtered probability space. 
%
Throughout the article, we use $\N$ to 
stand for the set of all positive integers.
%and denote $\N_0 = \{ 0 \} \cup \N$.
Let $ T \in (0, \infty) $ 
and let $C$ denote a generic deterministic positive constant, independent of the time stepsize $h  =  \tfrac{T}{N} >0$, $N \in \N$.
By $\| \cdot \|$ and $ \langle \cdot, \cdot \rangle $ we denote the Euclidean norm 
and the inner product of vectors in $\R^d, d \in \N$, 
respectively.  
%Adopting the same notation as the vector norm, we write $\|A\| : =\sqrt{\tr(A^{T}A)}$ as the trace norm of a matrix $A \in \R^{d \times m}$, $d,m \in \N$.
On a filtered probability space $ \left( \Omega, \mathcal{ F }, \{ \mathcal{ F }_t \}_{ t \in [0,T] }, \P 
\right) $, we use $\E$ to denote the expectation and $L^{r} (\Omega; \R^d ), r \geq 1 $, to denote 
the family of $\R^d$-valued random variables $\xi$ satisfying $\E[ \|\xi \|^{r}]<\infty$. Moreover, $W \colon [0,T] \times \Omega \rightarrow \mathbb{R}$ stands for the real-valued
standard Brownian motion with respect to $ \{ \mathcal{ F }_t \}_{ t \in [0,T] } $.

The present article aims to effectively compute
$
\mathbb{E} 
    \left[
    \phi\left( X ({T}) \right) 
    \right],
$
where $ \{ X(t) \}_{t \in [0, T]} $ is the unique solution of the following It\^o SDEs:
\begin{equation} \label{Milstein-eq:SODE}
\left\{
    \begin{array}{ll}
    \dd X (t)  = f ( X (t) ) \, \dd t + g ( X (t) ) \,\dd W (t),
    \quad
    t \in (0, T],
    \\
   X(0) = X_0 \in (0, \infty).
 \end{array}\right.
\end{equation}
Here the drift coefficient function $f  \colon \mathbb{R} \rightarrow \mathbb{R}$ and
the diffusion coefficient function $g \colon \mathbb{R} \rightarrow \mathbb{R}$
are, respectively, given by
\begin{equation} \label{eq:f-g-notation}
    f (x) := x ( \mu - \alpha x ),
    \quad
    g ( x ) : = \beta x^{3/2}.
\end{equation}
%assumed to obey a coupled monotonicity condition like \eqref{Milstein-eq:monotonicity-condition}.
%\begin{equation} \label{Milstein-eq:monotonicity-condition}
%\begin{split}
%\langle x - y, f ( t, x ) - f ( t, y ) \rangle 
%+ 
%\tfrac {q-1} {2} \lVert g ( t, x ) - g ( t, y ) \rVert ^ 2   \leq L \lVert x - y \rVert ^ 2,
%\quad
%\forall \, t\in [0, T], \: x, y \in \R^d,
%\end{split}
%\end{equation}
%for some $q > 2$.
%
%
Such an equation is referred to as the Heston $3/2$-model \cite{heston1997simple,ahn1999parametric}.
The well-posedness and regularity of the $3/2$-model
have been established in \cite{neuenkirch2014first,higham2013convergence}, which is recalled in the forthcoming proposition.
\begin{prop}
\label{Chpter5lem:32model-scheme-well-posedness}
Let $\mu , \alpha , \beta > 0$, $X_0 >0$.
Then the 
Heston $ 3/2 $-model \eqref{Milstein-eq:SODE}-\eqref{eq:f-g-notation} has a unique global positive solution.
Moreover, for $ 1 \leq  p \leq  1 + \tfrac{2 \alpha}{ \beta^2 }$
we have
\begin{equation}\label{eq:pro-wellposed-moment-bound}
\sup_{ t \in [0, T]} 
\E [ \| X (t) \|^p ]
\leq
C
( 1 + \| X_0 \|^p )
<
\infty.
\end{equation}
\end{prop}
To discretize \eqref{Milstein-eq:SODE}, we construct a uniform mesh $\{ t_i \}_{0 \leq i \leq N}$ 
on $[0, T]$ with $t_0 = 0$ and $ h = \tfrac{T}{N}$ being the stepsize, for any $ N \in \N$.
%
%%%%
%
On the uniform mesh, we propose the following split-step implicit Milstein type methods as follows:
\begin{equation}
\label{Chpter5Milstein-eq:scheme}
% \begin{split}
 \left\{
    \begin{array}{lll}
&    Z_n = Y_n  +  h f ( Z_n) -
            \tfrac {1} {2} h 
            g' (Z_n) g (Z_{n}),
    \\
&    Y_{ n + 1} =  Y_{ n } 
         + h f ( Z_n )
            + g(Z_{n}) \Delta W_n
            + 
            g' ( Z_n) g ( Z_{n} ) 
            I_{1,1}^{t_{n},t_{ n+1 }}
            ,
    \\
&   Y_{0}= X_0,
 \end{array}\right.
%\end{split}
\end{equation}
where, for $n = 0, 1, 2, \ldots, N - 1$,
%$\theta, \eta \in [0, 1]$,
\begin{equation}
\Delta W_{ n } : = W ( {t_{ n + 1 } }  ) - W ( {t_ {n} } ),
\quad
I_{1,1}^{t_{n},t_{ n+1 }} 
:= \tfrac12 [ ( \Delta W_n )^2 - h ].
\end{equation}
Obviously, the above scheme can be also written as
\begin{equation}
\label{Chpter5Milstein-eq:scheme2}
% \begin{split}
 \left\{
    \begin{array}{ll}
&    Z_n = Y_n  +  h f ( Z_n) -
            \tfrac {1} {2} h 
            g' (Z_n) g (Z_{n}),
    \\
&    Y_{ n + 1} =  Z_n
            + g(Z_{n}) \Delta W_n
            + 
            \tfrac12
            g' ( Z_n) g ( Z_{n} ) 
            ( \Delta W_n )^2,
    \\
&
    Y_0 = X_0       .
 \end{array}\right.
%\end{split}
\end{equation}
As shown below, the scheme can be solved explicitly and has
the positivity preserving property.
\begin{prop}
Let $\mu , \alpha , \beta > 0$, $X_0 >0$.
Given any initial value $Y_0 = X_0 > 0$
and for any stepsize $h = \tfrac{T}{N} > 0$, 
the scheme \eqref{Chpter5Milstein-eq:scheme}
or \eqref{Chpter5Milstein-eq:scheme2} admits unique, positive numerical solutions 
$\{ Y_n \}_{ n \in \N} $:
\begin{equation} \label{eq:Spli-step-doub-implicit-Mil}
% \begin{split}
 \left\{
    \begin{array}{ll}
&    Z_n
 = \frac{ - ( 1 - h \mu ) 
+ 
\sqrt{ ( 1 - h \mu )^2 +  h (3  \beta^2 + 4 \alpha) Y_n }
}
{ 2 h (\tfrac34  \beta^2 + \alpha) },
    \\
&    Y_{ n + 1} =  Z_n
            + \beta(Z_{n})^{\frac32} \Delta W_n
            + 
            \tfrac34
            \beta^2 ( Z_{n} )^2 
            ( \Delta W_n )^2  
            ,
    \quad Y_0 = X_0.
 \end{array}\right.
%\end{split}
\end{equation}
In other words, the scheme \eqref{Chpter5Milstein-eq:scheme} or \eqref{Chpter5Milstein-eq:scheme2} is positivity preserving unconditionally.
\end{prop}
{\bf Proof.}
In view of \eqref{eq:f-g-notation}, the iteration $Z_n$ 
follows a quadratic equation
\[
Z_n = Y_n + h Z_n ( \mu - \alpha Z_n )
-
\tfrac34 \beta^2 h ( Z_n )^2, 
\]
so for any $Y_n > 0$, one can straightforwardly and explicitly find a unique 
solution $Z_n > 0$ to the quadratic equation, given by the fist equation of \eqref{eq:Spli-step-doub-implicit-Mil}.  
By noting that $Z_n > 0$, after some elementary rearrangements, one can deduce from the second equation of \eqref{eq:Spli-step-doub-implicit-Mil} that
\begin{equation}
    Y_{n+1} = \tfrac{2}{3} Z_n 
    + 
    \Big(
    \tfrac{\sqrt{3}}{2} \beta Z_n \Delta W_n 
    +
    \tfrac{1}{\sqrt{3}} \sqrt{ Z_n }
    \Big)^2
    > 0
\end{equation}
for any step-size $h>0$.
This confirms that the scheme \eqref{Chpter5Milstein-eq:scheme} or \eqref{Chpter5Milstein-eq:scheme2} is positivity preserving unconditionally.
\hfill $\square$
%%%%%
%The monotonicity condition \eqref{Milstein-eq:monotonicity-condition} guarantees that the implicit methods ($\theta \eta  >0$)
%\eqref{Milstein-eq:scheme} admit unique $\mathcal{F}_{t_n} $-adapted solutions $Z_n, Y_n$  in $\R^d$ when 
%$h \in (0, \tfrac{1}{ L_0 })$ (see, e.g., \cite[Corollary 4.2]{beyn2016stochastic}). 
%(consult \cite[Theorem 3.1???]{beyn2016stochastic} for example). 
%
%One can consult \cite[Section 2]{beyn2016stochastic} and \cite[Proposition 5.4]{beyn2016stochastic} for the above assertions.
%
%%%%%%%%%

%\section{Order-one mean-square convergence 
%and applications to MLMC}
%
\section{MLMC and its optimal complexity}
\label{Milstein-sec:conv-rates-global-polynomial}
To approximate expectations of some functions of the solution to the Heston 3/2-model effectively,
we turn to multilevel Monte Carlo (MLMC) methods.
This section is thus devoted to the analysis of the 
complexity of MLMC methods. At first we revisit the
basic elements of MLMC.

\subsection{Basic elements of MLMC}
For simplicity of the presentation, we denote
\begin{equation} \label{def:functional-P-in-mlmc}
\Phi^X := \phi(X({T}))
\end{equation} 
and denote $\bar{\Phi}_{\ell}$ as the approximation of $\Phi^X$ using a numerical discretization with timestep 
\[
h_{\ell} = \frac{T}{2^{\ell}}, \quad \ell \geq 0.
\]
For some $L > 1$, the simulation of $\mathbb{E} [\bar \Phi_{L}]$ can be split into a series of levels of resolution as
\begin{equation} \label{equation:main-idea-of-mlmc}
\mathbb{E}
\big[
\bar{\Phi}_{L}
\big]
=
\mathbb{E}
\big[
\bar{\Phi}_{0}
\big]
+
\sum_{\ell=1}^{L}
\mathbb{E}
\big[
\bar{\Phi}_{\ell}-\bar{\Phi}_{\ell-1}
\big] .
\end{equation}

The key idea of MLMC is to independently estimate each of the expectations on the right-hand side of \eqref{equation:main-idea-of-mlmc}  in a way which minimises the overall variance for a given computational cost. To do so,
we let $\bar{Z}_{0}$ be an estimator of $\mathbb{E}[\bar{\Phi}_{0}]$ with $N_{0}$ Monte Carlo samples and $\bar{Z}_{\ell}$ be an estimator of $\mathbb{E}[\bar{\Phi}_{\ell}-\bar{\Phi}_{\ell-1}]$ with $N_{\ell}$ Monte Carlo samples, i.e.
\begin{equation} \label{definition:estimator-of-multilevel-estimator}
\bar{Z}_{\ell}
=\left\{
  \begin{array}{ll}
  N_{0}^{-1}\sum_{i=1}^{N_{0}}\bar{\Phi}_{0}^{(i)}, & \ell=0, \\
  N_{\ell}^{-1}\sum_{i=1}^{N_{\ell}}
  \left(
  \bar{\Phi}_{\ell}^{(i)}-\bar{\Phi}_{\ell -1}^{(i)} 
  \right), 
  & \ell>0.
  \end{array}\right.
\end{equation}
Therefore, the final multilevel estimator $\bar{Z}$ is given by the sum of the level estimators \eqref{definition:estimator-of-multilevel-estimator} as follows:
 \begin{equation} \label{equation:final-multilevel-estimator}
\bar{Z}
:=
\sum_{\ell=0}^{L} 
\bar{Z}_{ \ell }  .
\end{equation}
Since $\bar{\Phi}_{\ell}^{(i)}-\bar{\Phi}_{\ell -1}^{(i)} $  comes from two discrete approximations for the same sample, the diﬀerence would be small (due to strong convergence) on finer levels of resolution and so the variance is also small. As a result, very few samples are required on ﬁner levels and more on  coarser levels
to accurately estimate the expected value more efficiently.
%The key feature of MLMC is that $\bar{\Phi}_{\ell}^{i}-\bar{\Phi}_{\ell -1}^{i} $ is derived from the same numerical scheme, but is applied to different grids.
%------------------------------------
%As $\ell$ increases, the difference of $\bar{\Phi}_{\ell}^{i}-\bar{\Phi}_{\ell -1}^{i}$ decreases due to the strong convergence, so that very few samples are required on the finer level to make the estimation of \eqref{equation:main-idea-of-mlmc} more efficient.
%
%The key feature of MLMC is that $\Phi_{\ell}^{i}-\Phi_{\ell -1}^{i} $ comes from the same numerical scheme but under different grids.

In what follows, we revisit the MLMC complexity theorem due to Giles \cite{giles2008multilevel}.
%------------------------------------------------------MLMC theorem--------------------------------------------
\begin{thm} \label{theorem:mlmc-complexity-theorem}
Let $\Phi^X$, $\bar{Z}_{\ell}$ and  $\bar{Z}$ be defined as  \eqref{def:functional-P-in-mlmc}, \eqref{definition:estimator-of-multilevel-estimator} and \eqref{equation:final-multilevel-estimator}, respectively. Let $\bar \Phi_{\ell}$ be the corresponding level $\ell$ numerical approximation and let $\mathcal{C}_{\ell}$ denote the  computational complexity of $\bar{Z}_{\ell}$. Assume that there exist positive constants $\chi $, $\gamma$, $\theta$, $K_{1}$, $K_{2}$, $K_{3}$ such that $\chi \geq \frac{1}{2}\min \{\gamma , \theta\}$ and
\begin{equation} 
\begin{aligned}
&\text { (a) }
\quad
  \left|
    \mathbb{E}\big[\bar{\Phi}_{\ell}-\Phi^X\big]
  \right| 
  \leq K_{1} h_{\ell}^{\chi},  \\  \notag
&\text { (b) } 
\quad
\mathbb{E}\left[\bar{Z}_{\ell}\right]
=\left\{
  \begin{array}{ll}
  \mathbb{E}\big[\bar{\Phi}_{0}\big], & \ell=0 \\
  \mathbb{E}\big[\bar{\Phi}_{\ell}-\bar{\Phi}_{\ell-1}\big], & \ell>0
  \end{array}\right.  \\
&\text { (c) } 
\quad
  \text{Var} \left[\bar{Z}_{\ell}\right] \leq K_{2} N_{\ell}^{-1} h_{\ell}^{\gamma},  \\
&\text { (d) } 
\quad
  \mathcal{C}_{\ell} \leq K_{3} N_{\ell} h_{\ell}^{-\theta}. 
\end{aligned}
\end{equation}
Then  there exists a positive constant $K_{4}$ such that for any $0 < \epsilon < e^{-1}$ there are values $L$ and $N_{\ell}$ \hspace{0.05em} for which the multilevel estimator \eqref{equation:final-multilevel-estimator}
has a mean-square-error 
\begin{equation}
\text{MSE} := \mathbb{E}\left[(\bar{Z}-\mathbb{E}[\Phi^X])^{2}\right]
<\epsilon^{2} 
\end{equation}
with a  computational complexity $\mathcal{C}$ satisfying
\begin{equation} \label{eq:comlexity-thm-C}
\mathcal{C} 
\leq
\left\{
\begin{array}{ll}
K_{4} \epsilon^{-2}, & \gamma>\theta \\
K_{4} \epsilon^{-2}(\log (1 / \epsilon))^{2}, & \gamma=\theta \\ 
K_{4} \epsilon^{-2-(\theta-\gamma) / \chi}, & 0<\gamma<\theta.
\end{array}\right.
\end{equation}
\end{thm}

\subsection{Order-one mean-square convergence of the scheme}
The main focus of this subsection is to analyze the approximation error of the proposed scheme and recover the
expected mean-square convergence rate of order one.
In order to do this,  we introduce
an auxiliary process
$\widetilde{X} ( t_n )$, $ n = 0, 1, \ldots N$, implicitly
given by
\begin{equation}\label{eq:Chapter5Aux-X}
\widetilde{X} ( t_n ) = X ( t_n )
+
 h f (\widetilde{X} ( t_n ))
-
\tfrac {h} {2} 
            g' (\widetilde{X} ( t_n )) 
            g (\widetilde{X} ( t_n ))
.
\end{equation}
Equipped with the auxiliary process, we first derive
an upper bound of the mean-square approximation error.
Such an idea originally comes from \cite{wang2020mean}.
\begin{thm}
\label{Chaper5Milstein-thm:convergence-rate}
Let $\mu , \alpha , \beta > 0$, $X_0 >0$ 
and $\alpha > \tfrac32 \beta^2$. 
Assume that $\{ X({ t}) \}_{ t \in [0, T]}$ and $ \{Y_n\}_{0\leq n\leq N} $ are unique positive solutions to
the Heston $ 3/2 $-model \eqref{Milstein-eq:SODE}-\eqref{eq:f-g-notation} and the proposed scheme \eqref{Chpter5Milstein-eq:scheme}, respectively.
Then
\begin{equation} 
\begin{split}
\E [ \| X( {t_n } ) - Y_n  \| ^2 ]
\leq
C
\Big (
%\tfrac{q-1}{q - 2}
\sum_{ j = 1 }^{ n } \E [ \| \mathcal{R}_{ j } \| ^2 ]
+
%\tfrac{ (\theta^2 + 1)  }{ \theta^2 h}
h^{-1}
\sum_{ j = 1 }^{ n }
\E \big[ \big \| 
\E ( \mathcal{ R }^f_{j} | \mathcal{F}_{t_{j-1}} )
\big \| ^2 \big]
\Big )
%\leq
%C h^\frac12
,
\end{split}
\end{equation}
where $\mathcal{R}_{j}:=  \mathcal{R}^f_{j} + \mathcal{R}^g_{j}$, 
$j = 1, 2, \ldots, N$ 
with $\mathcal{R}^f_{j}$ and $\mathcal{R}^g_{j}$  given by
\begin{equation}
\begin{split}
\label{Chpter5Milstein-eq:Rj-defn}
\mathcal{R}^f_{j}
& :=
\int_{ t_{j-1} }^{ t_{j} }
\{
f ( X ( s ) ) - f ( \widetilde{X} ( t_{j-1} ) )
\}
\, \dd s,
\\
\mathcal{R}^g_{j}
& 
:= 
\int_{ t_{j-1} }^{ t_{j} }
\{
g ( X ( s ) ) - g ( \widetilde{X} ( t_{j-1} ) )
\}
\, \dd  W ( s )
-
            g'g (  \widetilde{X} ( t_{j-1} ) ) I_{1,1}^{t_{j-1},t_{j}}
.
\end{split}
\end{equation}
\end{thm}
To prove the theorem, we rely on the following lemma.
%%%%%
\begin{lem}
\label{Milstein-lem:monotonicity-condition}
Let $\mu , \alpha , \beta > 0$, $X_0 >0$ 
and $\alpha > \tfrac32 \beta^2$. 
Then there exist non-negative constants $ 2 < q < 1 + \tfrac{ 8 \alpha }{ 9 \beta^2 }  $, $ \varrho \in (1, \infty) $, $L_0, L_1 \in [0, \infty)$, $ h_0 \in (0, T]$ such that, for any $x, y \in (0, \infty)$,
\begin{equation} \label{Milstein-eq:monotonicity-condition}
\begin{split}
& \big \langle x - y , [ f ( x ) - f ( y ) ]
 - \tfrac{1}{2} 
  [ g' (x) g( x )
  - g' (y) g( y )]  \big \rangle
   \leq L_0 \| x - y \|^2,
\\
& 2 \langle x - y ,f ( x ) - f ( y ) \rangle
  +  ( q - 1 ) \|g ( x ) - g ( y )\| ^2
  + \tfrac{ \varrho h }{ 2 }  
  \big\| 
   g'( x ) g ( x ) - g'( y ) g ( y ) \big\| ^2
 \nonumber 
 \\ \:
 & \quad
  +  h  \big \langle 
  \big[ g'( x ) g ( x )
   - g'( y ) g ( y ) \big] , f ( x ) - f ( y ) \big \rangle
  - h \| f ( x ) - f ( y ) \|^2
%  \nonumber \\ 
% & \quad 
   \leq L _ 1 \| x - y \|^2.
\end{split}
\end{equation}
\end{lem}
{\bf Proof.}
For short we denote
\begin{equation}
\Theta (x, y, h)
:=
\tfrac{ \varrho h }{ 2 } \| g'g ( x ) - g'g ( y ) \| ^2
  + h \langle  g'g ( x ) - g'g ( y ) , f ( x ) - f ( y ) \rangle
  - h \| f ( x ) - f ( y ) \|^2.
\end{equation} 
Noting  $ g ' (x) g ( x ) = \tfrac32 \beta^2 x^2, \, x \in \R_+ $, one can find a constant 
$\tilde{c} > 0$ such that
\begin{equation}
\begin{split}
\Theta (x, y, h)
&  =
\tfrac98 \varrho \beta^4 h ( x^2 - y^2 )^2 + \tfrac32 \beta^2 \mu h ( x^2 - y^2 ) ( x - y)
 - \tfrac32 \beta^2 \alpha h ( x^2 - y^2 )^2
 \\
& \quad
 - h \big[ \mu^2 ( x - y)^2
 - 2 \mu \alpha ( x - y ) ( x^2 - y^2 ) + \alpha^2 ( x^2 - y^2 )^2 \big]
\\
&  = 
\big ( [ \tfrac98 \varrho \beta^4 - \tfrac32 \beta^2 \alpha - \alpha^2 ] ( x + y )^2 
+ 
[ \tfrac32 \beta^2 \mu +  2 \mu \alpha ] ( x + y ) - \mu^2   \big )  ( x - y)^2 h
\\
& \leq
%\big[ \tfrac32 \beta^2 \mu +  2 \mu \alpha \big] ( x + y )  ( x - y)^2 h,
\tilde{c} ( x - y)^2 h,
\qquad
\forall \, x , y \in \R_+,
\end{split}
\end{equation}
where we used the fact that $ \tfrac98 \varrho \beta^4 - \tfrac32 \beta^2 \alpha - \alpha^2 < 0$ for some $\varrho > 1$, by the assumption
$ 
\alpha > \tfrac32 \beta^2.
%> \tfrac{3(\sqrt{3} -1)}{4} \beta^2
$
Moreover, we take $2 < q < 1 + \tfrac{ 8 \alpha }{ 9 \beta^2 } $ to promise $ \tfrac94 (q - 1 ) \beta^2 - 2 \alpha \leq 0$.
As a result,
\begin{equation}
\begin{split}
&
2 \langle x - y ,f ( x ) - f ( y ) \rangle
  +  ( q - 1 ) \|g ( x ) - g ( y )\| ^2
+
\Theta (x, y, h)
\\
& = 
2 \mu | x - y |^2 -  2\alpha (x^2 - y^2) ( x - y) 
+
(q - 1 ) \beta^2 ( x^\frac32 - y^\frac32 )^2
+
\Theta (x, y, h)
\\
& \leq
2\mu | x - y |^2 
+
% \big[ \tfrac32 \beta^2 \mu h +  2 \mu \alpha h 
% - 2 \alpha + 3 (q - 1 ) \beta^2 \big] ( x + y )  ( x - y)^2
\big[ 
\tfrac94 (q - 1 ) \beta^2  - 2 \alpha \big] ( x + y )  ( x - y)^2
+
\tilde{c} ( x - y)^2 h
\\
& \leq
(2 \mu + \tilde{c} T ) | x - y |^2,
\quad
 \forall \, x , y \in \R_+,
\end{split}
\end{equation}
which validates the desired assertion.
\hfill $\square$

At the moment, we are well prepared to prove 
Theorem \ref{Chaper5Milstein-thm:convergence-rate}

\textbf{Proof of Theorem \ref{Chaper5Milstein-thm:convergence-rate}.}
By \eqref{Chpter5Milstein-eq:Rj-defn}, one infers
\begin{equation}
\begin{split}
\label{Chaper5Milstein-eq:Xt-Rn}
X ( t_{ n + 1 } )
& = X ( t_{ n } )
+
h f ( \widetilde{X} ( t_n ))
+
g ( \widetilde{X} ( t_n )) \Delta W_{ n }
+
            g' g ( \widetilde{X} ( t_{n} ) ) I_{1,1}^{t_{n},t_{ n+1 }}
+
\mathcal{R}_{n+1}.
\end{split}
\end{equation}
To ease the notation, we denote
\begin{equation}
\begin{split}
e_n  &: = X( {t_n } ) - Y_n,
\\
\delta f_{ \widetilde{X} ( t_n ), Z_n }
&: =
f ( \widetilde{X} ( t_n ))
-
f ( Z_n),
\\
\delta g_{ \widetilde{X} ( t_n ), Z_n }
& : =
g ( \widetilde{X} ( t_n ))
-
g ( Z_n),
\\
\delta ( g'g )_{ \widetilde{X} ( t_{n} ), Z_n }
& :=
g'g ( \widetilde{X} ( t_n ))
-
g'g ( Z_n).
\end{split}
\end{equation}
Subtracting the first equation of \eqref{Chpter5Milstein-eq:scheme} from \eqref{eq:Chapter5Aux-X} yields
\begin{equation}
\label{Chapter5Milstein-eq:e_n-delta-f}
\begin{split}
\widetilde{X} ( t_n ) - Z_n
& =
e_n
+
h [
f ( \widetilde{X} ( t_n ))
-
f ( Z_n)
]
-
\tfrac{1}{2} h
[ g'g ( \widetilde{X} ( t_n ) ) - g' g ( Z_n ) ]
\\
&
=
e_n
+
h
\delta f_{ \widetilde{X} ( t_n ), Z_n }
-
\tfrac{1}{2} h
\delta (g'g ) _{ \widetilde{X} ( t_n ), Z_n }
.
\end{split}
\end{equation}
In addition, 
subtracting the second equation of \eqref{Chpter5Milstein-eq:scheme} from
\eqref{Chaper5Milstein-eq:Xt-Rn} and using the above notation result in
\begin{equation} \label{Milstein-eq:error-equation}
\begin{split}
e_{n+1}
%&
%=
%e_n
%+
%h
%[
%f ( t_n + \theta h,  \widetilde{X} ( t_n ))
%-
%f ( t_n + \theta h,  Z_n)
%]
%\\ & \quad
%+
%[
%g ( t_n + \theta h,  \widetilde{X} ( t_n ))
%-
%g ( t_n + \theta h,  Z_n)
%]
%\Delta W_{ n }
%+
%\mathcal{R}_{ n + 1 }
%\\
%&
=
e_n
+
h
\delta f_{ \widetilde{X} ( t_n ), Z_n }
+
\delta g_{ \widetilde{X} ( t_n ), Z_n }
\Delta W_{ n }
+
            \delta ( g' g )_{ \widetilde{X} ( t_{n} ), Z_n } I_{1,1}^{t_{n},t_{ n+1 }}
+
\mathcal{R}_{ n + 1 }.
\end{split}
\end{equation}
As a result,
\begin{equation} \label{Chapter5Milstein-eq:e-n2-identity}
\begin{split}
\| e_{ n + 1} \|^2
= &
\| e_{ n } \|^2
+
h^2 \| \delta f_{ \widetilde{X} ( t_n ), Z_n } \|^2
+
\| \delta g_{ \widetilde{X} ( t_n ), Z_n } \Delta W_n \|^2
\\
& \quad
+
\Big \|
        \delta ( g'g )_{ \widetilde{X} ( t_{n} ), Z_n } I_{1,1}^{t_{n},t_{ n+1 }}
\Big\|^2
+
\| \mathcal{ R }_{ n+ 1} \|^2
\\
& \quad
+
2 h \langle e_n, \delta f_{ \widetilde{X} ( t_n ), Z_n } \rangle
+
2 \langle e_n, \delta g_{ \widetilde{X} ( t_n ), Z_n } \Delta W_n \rangle
\\
& \quad
+
2 \langle e_n, 
            \delta ( g'g )_{ \widetilde{X} ( t_{n} ), Z_n } I_{1,1}^{t_{n},t_{ n+1 }} \rangle
+
2 \langle e_n, \mathcal{ R }_{n+1} \rangle
\\
& \quad
+
2 h \langle \delta f_{ \widetilde{X} ( t_n ), Z_n },  \delta g_{ \widetilde{X} ( t_n ), Z_n } \Delta W_n \rangle
\\
& \quad
+
2 h \langle \delta f_{ \widetilde{X} ( t_n ), Z_n },
            \delta ( g'g )_{ \widetilde{X} ( t_{n} ), Z_n } I_{1,1}^{t_{n},t_{ n+1 }} \rangle
\\
& \quad
+
2 h \langle \delta f_{ \widetilde{X} ( t_n ), Z_n },  \mathcal{ R }_{n+1} \rangle
\\
& \quad
+
2 \langle \delta g_{ \widetilde{X} ( t_n ), Z_n } \Delta W_n,
            \delta ( g'g )_{ \widetilde{X} ( t_{n} ), Z_n } I_{1,1}^{t_{n},t_{ n+1 }}
\rangle
\\
& \quad
+
2 \langle \delta g_{ \widetilde{X} ( t_n ), Z_n } \Delta W_n, \mathcal{ R }_{n+1} \rangle
+
2 \big \langle
            \delta ( g'g )_{ \widetilde{X} ( t_{n} ), Z_n } I_{1,1}^{t_{n},t_{ n+1 }},
            \mathcal{ R }_{n+1} \big \rangle
.
\end{split}
\end{equation}
By observing that
\begin{equation}
\begin{split}
& \E [ \langle e_n, \delta g_{ \widetilde{X} ( t_n ), Z_n } \Delta W_n \rangle ] = 0,
\quad
\E [ \langle \delta f_{ \widetilde{X} ( t_n ), Z_n },  \delta g_{ \widetilde{X} ( t_n ), Z_n } \Delta W_n \rangle ] = 0,
\\
& \E \Big[  \Big \langle e_n, 
            \delta ( g'g )_{ \widetilde{X} ( t_{n} ), Z_n } I_{1,1}^{t_{n},t_{ n+1 }} \Big \rangle \Big]
= 0,
\\
&
\E \Big [
   \Big \langle \delta f_{ \widetilde{X} ( t_n ), Z_n }, 
            \delta ( g'g )_{ \widetilde{X} ( t_{n} ), Z_n } I_{1,1}^{t_{n},t_{ n+1 }}
            \Big \rangle
      \Big]
= 0,
\\
& \E \Big[
\Big \langle \delta g_{ \widetilde{X} ( t_n ), Z_n } \Delta W_n,
            \delta ( g'g )_{ \widetilde{X} ( t_{n} ), Z_n } I_{1,1}^{t_{n},t_{ n+1 }}
\Big \rangle
\Big]
=
0,
\end{split}
\end{equation}
we further deduce
\begin{equation} \label{eq:error-recurrence}
\begin{split}
\E [ \| e_{ n + 1} \|^2 ]
& =
\E [ \| e_{ n } \|^2 ]
+
h^2  \E [ \| \delta f_{ \widetilde{X} ( t_n ), Z_n } \|^2 ]
+
h \E [ \| \delta g_{ \widetilde{X} ( t_n ), Z_n } \|^2 ]
\\
& \quad
+
\tfrac{h^2}{2}
\E
\big[ \big \|
            \delta ( g'g )_{ \widetilde{X} ( t_{n} ), Z_n }
\big\|^2
\big]
+
\E [ \| \mathcal{ R }_{ n+ 1} \|^2 ]
+
2 h \E [ \langle e_n, \delta f_{ \widetilde{X} ( t_n ), Z_n } \rangle ]
\\
& \quad
+
2
\E \big[ \langle e_n, \mathcal{ R }_{n+1} \rangle \big]
+
2 h \E [ \langle \delta f_{ \widetilde{X} ( t_n ), Z_n },  \mathcal{ R }_{n + 1}  \rangle ]
\\
& \quad
+
2 \E [
\langle \delta g_{ \widetilde{X} ( t_n ), Z_n } \Delta W_n, \mathcal{ R }_{n+1} \rangle
]
\\
& \quad
+
2 \E \big[
\big\langle 
            \delta ( g'g )_{ \widetilde{X} ( t_{n} ), Z_n } I_{1,1}^{t_{n},t_{ n+1 }},
            \mathcal{ R }_{n+1} \big \rangle \big]
.
\end{split}
\end{equation}
With the aid of \eqref{Chapter5Milstein-eq:e_n-delta-f}, we get
\begin{equation}
\begin{split}
2 h \langle e_n, \delta f_{ \widetilde{X} ( t_n ), Z_n } \rangle
& =
2 h
\langle \widetilde{X} ( t_n ) - Z_n
-
h
\delta f_{ \widetilde{X} ( t_n ), Z_n }
%\\
%& \quad
+
\tfrac{1}{2} h
\delta (g'g ) _{ \widetilde{X} ( t_n ), Z_n }
,
\delta f_{ \widetilde{X} ( t_n ), Z_n } \rangle
\\
& =
2 h
\langle \widetilde{X} ( t_n ) - Z_n,
\delta f_{ \widetilde{X} ( t_n ), Z_n } \rangle
- 2 h^2 \| \delta f_{ \widetilde{X} ( t_n ), Z_n } \|^2
\\
& \quad
+
h^2 
\langle
\delta (g'g ) _{ \widetilde{X} ( t_n ), Z_n },
\delta f_{ \widetilde{X} ( t_n ), Z_n }
\rangle
,
\\
2 h  \langle \delta f_{ \widetilde{X} ( t_n ), Z_n },  \mathcal{ R }_{n+1} \rangle
& =
2
\big \langle
\widetilde{X} ( t_n ) - Z_n - e_n
+
\tfrac{1}{2} h
\delta (g'g ) _{ \widetilde{X} ( t_n ), Z_n }
,
\mathcal{ R }_{n+1}
\big \rangle
%\\
%& =
%...
.
\end{split}
\end{equation}
Plugging these two identities into \eqref{eq:error-recurrence} gives
\begin{equation}
\label{Chpter5Milstein-eq:en1_en_relation}
\begin{split}
\E [ \| e_{ n + 1} \|^2 ]
& =
\E [ \| e_{ n } \|^2 ]
- h^2  \E [ \| \delta f_{ \widetilde{X} ( t_n ), Z_n } \|^2 ]
+
2 h \E [ \langle \widetilde{X} ( t_n ) - Z_n, \delta f_{ \widetilde{X} ( t_n ), Z_n } \rangle ]
\\
& \quad
+
h \E [ \| \delta g_{ \widetilde{X} ( t_n ), Z_n } \|^2 ]
+
h^2 
\E \big[\langle
\delta (g'g ) _{ \widetilde{X} ( t_n ), Z_n },
\delta f_{ \widetilde{X} ( t_n ), Z_n }
\rangle 
\big]
\\
& \quad
+
\tfrac{h^2}{2}
\E
\big[ \big \|
            \delta ( g'g )_{ \widetilde{X} ( t_{n} ), Z_n }
\big\|^2
\big]
+
\E [ \| \mathcal{ R }_{ n+ 1} \|^2 ]
\\
& \quad
+
2 \E [ \langle \widetilde{X} ( t_n ) - Z_n,  \E ( \mathcal{ R }_{n+1} | \mathcal{F}_{t_n} ) \rangle ]
+
h
\E [ \langle \delta ( g'g )_{ \widetilde{X} ( t_n ), Z_n },   \mathcal{ R }_{n+1} \rangle ]
\\
& \quad
+
2 \E [
\langle \delta g_{ \widetilde{X} ( t_n ), Z_n } \Delta W_n, \mathcal{ R }_{n+1} \rangle
]
\\ & \quad
+
2 \E \big[
\big\langle
            \delta ( g'g )_{ \widetilde{X} ( t_{n} ), Z_n } I_{1,1}^{t_{n},t_{ n+1 }},
            \mathcal{ R }_{n+1} \big \rangle \big].
\end{split}
\end{equation}
Owing to the Young inequality and noticing that 
\begin{equation}
\E ( \mathcal{ R }_{n+1} | \mathcal{F}_{t_n} ) = \E ( \mathcal{ R }^f_{n+1} | \mathcal{F}_{t_n} ),
\end{equation}
we have
\begin{equation}
\begin{split}
\E [ \| e_{ n + 1} \|^2 ]
& \leq
\E [ \| e_{ n } \|^2 ]
- h^2  \E [ \| \delta f_{ \widetilde{X} ( t_n ), Z_n } \|^2 ]
+
2 h \E [ \langle \widetilde{X} ( t_n ) - Z_n, \delta f_{ \widetilde{X} ( t_n ), Z_n } \rangle ]
\\
& \quad
+
h ( q - 1) \E [ \| \delta g_{ \widetilde{X} ( t_n ), Z_n } \|^2 ]
+
h^2 
\E \big[\langle
\delta (g'g ) _{ \widetilde{X} ( t_n ), Z_n },
\delta f_{ \widetilde{X} ( t_n ), Z_n }
\rangle 
\big]
\\
& \quad
+
\tfrac{ \varrho h^2}{2}
\E
\big[ \big \|
            \delta ( g'g )_{ \widetilde{X} ( t_{n} ), Z_n }
\big\|^2
\big]
+
\E [ \| \mathcal{ R }_{ n+ 1} \|^2 ]
\\
& \quad
+
h \E [ \| \widetilde{X} ( t_n ) - Z_n \|^2 ]
+
h^{-1} \E [ \| \E ( \mathcal{ R }^f_{n+1} | \mathcal{F}_{t_n} )  \|^2 ]
\\
& \quad
+
\tfrac{1} {q-2}\E [
\| \mathcal{ R }_{n+1} \|^2
]
+
\tfrac{1}{ \varrho - 1} \E \big[ \|
            \mathcal{ R }_{n+1} \|^2 \big].
\end{split}
\end{equation}
Using the second monotoncity condition in Lemma \ref{Milstein-lem:monotonicity-condition} ensures
\begin{equation} \label{eq:en-last-but-one}
\begin{split}
    \E [ \| e_{ n + 1} \|^2 ]
& \leq
\E [ \| e_{ n } \|^2 ]
+ h ( L_1 + 1 ) \E [ \| \widetilde{X} ( t_n ) - Z_n  \| ^2 ]
\\
& \quad
+
\big(
\tfrac{q-1}{q-2}
+
\tfrac{1}{ \varrho - 1} 
\big)
\E \big[ \|
            \mathcal{ R }_{n+1} \|^2 \big]
+
h^{-1} \E [ \|  \E ( \mathcal{ R }^f_{n+1} | \mathcal{F}_{t_n} ) \|^2 ].
\end{split}
\end{equation}
Further, one can use the first monotoncity condition in Lemma \ref{Milstein-lem:monotonicity-condition} to 
derive from \eqref{Chapter5Milstein-eq:e_n-delta-f} that
\begin{equation}
\begin{split}
\| e_n \|^2
 & =
\| \widetilde{X} ( t_n ) - Z_n - h
\delta f_{ \widetilde{X} ( t_n ), Z_n } 
+
\tfrac{1}{2} h
\delta (g'g ) _{ \widetilde{X} ( t_n ), Z_n }
\|^2
\\
&
\geq
\| \widetilde{X} ( t_n ) - Z_n \|^2
- 2 
\langle \widetilde{X} ( t_n ) - Z_n,  h \delta f_{ \widetilde{X} ( t_n ), Z_n } - \tfrac{1}{2} h
\delta (g'g ) _{ \widetilde{X} ( t_n ), Z_n } \rangle
 \\ &
\geq
( 1 - 2  L_0 h ) \| \widetilde{X} ( t_n ) - Z_n \|^2,
\end{split}
\end{equation}
which, on the condition $1 - 2 L_0 h > 0$, implies
\begin{equation}
\label{Chaper5Milstein-eq:tildeX-Z-en}
\| \widetilde{X} ( t_n ) - Z_n \|^2
\leq
\tfrac{1}{1 - 2 L_0 h}
\| e_n \|^2
.
\end{equation}
This together with \eqref{eq:en-last-but-one} shows
\begin{equation}
\begin{split}
    \E [ \| e_{ n + 1} \|^2 ]
& \leq
\big(1 + \tfrac{ h ( L_1 + 1 )}{1 - 2 L_0 h}  \big)
\E [ \| e_{ n } \|^2 ]
\\
& \quad
+
\big(
\tfrac{q-1}{q-2}
+
\tfrac{1}{ \varrho - 1} 
\big)
\E \big[ \|
            \mathcal{ R }_{n+1} \|^2 \big]
+
h^{-1} \E [ \|  \E ( \mathcal{ R }^f_{n+1} | \mathcal{F}_{t_n} ) \|^2 ].
\end{split}
\end{equation}
The proof is finished by iteration. 
\hfill $\square$
%

%%%%%%%%%%%%%%%%%%%%%%%%%
In light of \eqref{eq:f-g-notation}, one can easily compute that
\begin{align}
\label{Milstein-eq:drift-polynomial-growth-condition2}
\| f ( x ) - f ( y ) \| & \leq 
( \mu + \alpha \| x \|  + \alpha \| y \| ) \| x - y \|,
\quad
\forall  \, x, y >0,
\\ 
\| g ( x ) - g ( y ) \|^2 
&
\leq
\tfrac94 \beta^2
( \| x\| + \| y \| ) \| x - y \|^2,
\quad
\forall \, x, y > 0,
\\
\| g'(x) g ( x ) - g' (y) g ( y ) \|^2 
&
\leq
\tfrac94 \beta^4
\|
 x
+
y
\|^2
\|
 x
-
y
\|^2,
\quad
\forall \, x, y > 0.
\end{align}
and that
\begin{equation} \label{Milstein-eq:f-g-polynomial-growth}
\| f ( x ) \| \leq
 \alpha \| x \|^2,
\quad
\| g ( x) \| \leq
\beta \| x \|^{ \frac{ 3 } { 2 } } 
,
\quad
\| g'(x) g ( x) \| \leq
\tfrac32 \beta^2 \| x \|^2
,
\quad
\forall \, x >0 .
\end{equation}
Moreover, for  $p \leq 1 + \tfrac{2 \alpha}{\beta^2} $,
\begin{align} 
\label{Milstein-assum:eq-monotone-growth}
\langle x , f ( x )  \rangle + \tfrac{ p - 1}{2} \lVert 
g ( x ) \rVert ^ 2 
 =
\mu x^2 + ( \tfrac{p-1}{2} \beta^2 - \alpha ) x^3
\leq  \mu \|x\|^2,
\quad
\forall \, x >0.
%\end{equation}
\end{align}
These facts together with \eqref{eq:pro-wellposed-moment-bound} enable us to derive for any $ t , s \in [ 0, T ]$
and $p \leq 1 + \tfrac{2 \alpha}{\beta^2} $,
\begin{equation}\label{Milstein-eq:Xt1-Xt2}
\left \lVert X ( { t } ) - X ( { s } ) \right \rVert_{L^\delta ( \Omega; \R^d )}
\leq 
C  \lvert t - s  \rvert ^{ \frac 12 },
\quad
\delta \in [1, \tfrac { p } { 2 }  ].
\end{equation}
%%%%%%%%%%%%%%%%%%%%%%%%%%
%
In what follows, we obtain the order-one mean-square convergence of the proposed time discretization scheme.
\begin{thm}
\label{thm:Chapter5Milstein-convergence-rate-32model}
Let $X_0 > 0$ and $\mu , \alpha , \beta > 0$ 
satisfying $ \alpha \geq \tfrac52 \beta^2 $. 
Assume that $\{ X({ t}) \}_{ t \in [0, T]}$ and $ \{Y_n\}_{0\leq n\leq N} $ are unique positive solutions to
the Heston $ 3/2 $-model \eqref{Milstein-eq:SODE}-\eqref{eq:f-g-notation} and the proposed scheme \eqref{Chpter5Milstein-eq:scheme}, respectively.
For $ h \in (0, \tfrac{ \vartheta }{\mu} ] $ with
some $\vartheta \in (0, 1)$,
there exists a constant $C>0$ independent if 
$N \in \N$ such that
\begin{equation} \label{eq:Chapter5thm-32model-converg-rate}
	\begin{split}
	\sup_{ 0 \leq n \leq N} \E [ \lVert Y_n - X(t_n) \rVert^2 ]
	\leq &
	C h^2
%	\big( 1 +  \| X_0 \|_{L^6 ( \Omega; \R) }^3 \big)
	.
	\end{split}
	\end{equation}
\end{thm}
%
%%%
{\bf Proof.}
Following the above notation \eqref{eq:f-g-notation},
we get
\begin{equation}
\| x - h f(x) + \tfrac h 2 g' (x) g (x) \|^2
=
\| x - h x ( \mu - \alpha x  )
+
\tfrac{3h}{4} \beta^2 x^2
\|^2
\geq
( 1 - \mu h)^2 x^2,
\quad
\forall x >0.
\end{equation}
%%%%%%%%%%%%%%%%%%%%%%
%\begin{equation}
%\big \langle x, f(x) - \tfrac12 g' (x) g (x) \big \rangle
%=
%x^2 ( \mu - \alpha x  )
%-
%\tfrac34 \beta^2 x^3
%\leq
%\mu x^2,
%\quad
%\forall x >0.
%\end{equation}
Evidently, $  \widetilde{X} ( 0 ) = Z_0$ and $\widetilde{X} ( t_n )$  defined by \eqref{eq:Chapter5Aux-X} is positive as
$X ( t_n ) > 0$.
Then
\begin{equation} \label{Chapter5Milstein-eq:X-tilde-X-estimate}
\begin{split}
\| X ( t_n ) \|^2 
 & = 
\big\|  
\widetilde{X} ( t_n ) - h f (\widetilde{X} ( t_n )) 
+
\tfrac {h} {2} 
            g' (\widetilde{X} ( t_n )) 
            g (\widetilde{X} ( t_n ))
\big\|^2 
\\ &
\geq
(1 - \mu h )^2 \|  \widetilde{X} ( t_n ) \|^2
.
\end{split}
\end{equation}
%%%
%\begin{lem} \label{Milstein-lem:X-tilde-X-relationship}
%Suppose that Assumptions \ref{Milstein-ass:monotonicity-condition}, \ref{Milstein-ass:f-polynomial-growth} hold and 
For some $0 < \vartheta <1 $,
$ \mu h \leq  \vartheta $, one knows
\begin{equation} 
\label{Chapter5Milstein-eq:widetilde-X-X-relation}
\|  \widetilde{X} ( t_n ) \|
\leq
\tfrac{ 1 }{  1 - \mu h }
\| X ( t_n ) \|
\leq
\tfrac{ 1 }{  1 - \vartheta }
\| X ( t_n ) \|.
\end{equation}
Using the H\"{o}lder inequality shows
\begin{equation} \label{eq:Rj2-decomp-32model}
\begin{split}
\E [ \| \mathcal{R}_{ j } \|^2 ]
& \leq
2
h
\int_{ t_{j-1} }^{ t_{j} }
\E [
\|
f ( X ( s ) ) - f ( \widetilde{X} ( t_{j-1} ) ) \|^2
]
\, \dd s
\\
& \quad +
2
\E
\bigg [ \bigg\|
\int_{ t_{j-1} }^{ t_{j} }
\{
g ( X ( s ) ) - g ( \widetilde{X} ( t_{j-1} ) )
\}
\, \dd  W ( s )
%\\ & \qquad  \qquad
 - 
 g'g (  \widetilde{X} ( t_{j-1} ) ) I_{1,1}^{t_{j-1},t_{j}}
\bigg \|^2 \bigg ]
\\
& =: \mathbb{J}_1 + \mathbb{J}_2.
\end{split}
\end{equation}
%
%Similar to the estimate of $J_{12}$,  $J_{13}$  in \eqref{Milstein-eq:J12-estimate} and \eqref{Milstein-eq:J13-estimate}, one can easily bound $\mathbb{J}_1$ as follows:
Owing to \eqref{Milstein-eq:drift-polynomial-growth-condition2} and \eqref{eq:Chapter5Aux-X},
one can employ the H\"{o}lder inequality to arrive at
\begin{equation}
\label{Chaper5Milstein-eq:Err1}
\begin{split}
\mathbb{J}_1
& \leq
4 h
\int_{ t_{j-1} }^{ t_{j} }
\E [
\|
f (  X ( s ) ) - f ( X ( { t_{j-1} } ) )  \|^2
]
\, \dd s
\\ & \quad
+
4 h
\int_{ t_{j-1} }^{ t_{j} }
\E [
\|
f (  X ( { t_{j-1} } ) ) - f ( \widetilde{X} ( t_{j-1} ) )  \|^2
]
\, \dd s
\\ &
\leq
4 h
\int_{ t_{j-1} }^{ t_{j} }
\E [
( \mu + \alpha \|  X ( s ) \| +  \alpha \| X ( { t_{j-1} } ) \| )^2 \| X ( s ) - X ( { t_{j-1} } ) \| ^2
]
\, \dd s
\\ &
\quad
+
4 h
\int_{ t_{j-1} }^{ t_{j} }
\E [
( \mu + \alpha \|  \widetilde{X} ( t_{j-1} ) \| +  \alpha \| X ( { t_{j-1} } ) \| )^2 \| \widetilde{X} ( t_{j-1} ) - X ( { t_{j-1} } ) \| ^2
]
\, \dd s
\\ &
\leq
C
\Big(
1 +  \sup_{s \in [0, T] } \| X ( s ) \|_{ L^{ 6 } ( \Omega; \R^d ) }^{ 6 }
\Big) h^3.
\end{split}
\end{equation}
Thanks to the It\^{o} isometry,
\begin{equation}
\begin{split}
\mathbb{J}_2
&=
    2
\E
\bigg [ \bigg\|
\int_{ t_{j-1} }^{ t_{j} }
\{
g ( X ( s ) ) - g ( \widetilde{X} ( t_{j-1} ) )
- 
 g'g (  \widetilde{X} ( t_{j-1} ) ) ( W (s) - W ( t_{j-1} ) )
\}
\, \dd  W ( s )
\bigg \|^2 \bigg ]
\\
&
=
2
\int_{ t_{j-1} }^{ t_{j} }
\E
[
\|
g ( X ( s ) ) - g ( \widetilde{X} ( t_{j-1} ) )
- 
 g'g (  \widetilde{X} ( t_{j-1} ) ) ( W (s) - W ( t_{j-1} ) )
\|^2
]
\, \dd  s
\\ &
\leq
6
\int_{ t_{j-1} }^{ t_{j} }
\E
[
\|
g ( X ( s ) ) - g ( X ( t_{j-1} ) )
- 
 g'g (  X ( t_{j-1} ) ) ( W (s) - W ( t_{j-1} ) )
\|^2
]
\, \dd  s
\\ &
\quad
+
6 \int_{ t_{j-1} }^{ t_{j} }
\E
[
\|
g ( \widetilde{X} ( t_{j-1} ) )
- 
g ( X ( t_{j-1} ) )
\|^2
]
\, \dd  s
\\ &
\quad
+
6
\int_{ t_{j-1} }^{ t_{j} }
\E
[
\|
 [ g'g (  \widetilde{X} ( t_{j-1} ) ) 
- g'g (  X ( t_{j-1} ) ) ] ( W (s) - W ( t_{j-1} ) )
\|^2
]
\, \dd  s
\\
&
=:
\mathbb{J}_{21}
+
\mathbb{J}_{22}
+
\mathbb{J}_{23}.
\end{split}
\end{equation}
For $g (x) = \beta x^{ \frac32 }$ and 
$g' ( x )g (x) = \tfrac32 \beta^2 x^2$,  $x \in \R_+$, we employ the It\^o isometry and the It\^o formulae to obtain
%, 
%when $ \alpha \geq \tfrac52 \beta^2 $,
%%
\begin{equation} 
\begin{split}
\mathbb{J}_{21}
& =
6
\int_{ t_{j - 1 } }^{ t_{ j } } 
\E \Big[
\Big \lVert
\int_{ t_{j - 1 } }^{ s }
[
g' ( X(r) ) f ( X(r) ) 
+ \tfrac12 g'' ( X (r) ) g^2 ( X(r) ) ] \, \dd r
\\ &
\qquad \qquad
+
\int_{ t_{j - 1 } }^{ s }
[
g' g ( X(r) ) - g'g ( X({ t_{j-1} }) )
] \, \dd W_r
\Big \rVert^2 \Big]
\dd s
\\
&
\leq
12 h
\int_{ t_{j - 1 } }^{ t_{ j } } 
\int_{ t_{j - 1 } }^{ s }
\E
[
\|
g' ( X (r) ) f ( X (r) ) 
+ \tfrac12 g'' ( X_r ) g^2 ( X_r ) 
\|^2
]
\, \dd r \dd s
\\ & \quad
+
12
\int_{ t_{j - 1 } }^{ t_{ j } } 
\int_{ t_{j - 1 } }^{ s }
\E
[
\|
g' g ( X (r) ) - g'g ( X({ t_{j-1} }) )
\|^2
]
\, \dd r \dd s
\\ 
&
\leq
C h^3 \big( 1 + \sup_{s \in [0, T]} \| X_s \|_{L^6 ( \Omega; \R) }^6 \big).
%\\ &
%\leq
%C h^3.
\end{split}
\end{equation}
Here $ g' (x) f (x) = \tfrac32 \beta x^\frac32 ( \mu - \alpha x  ) $, 
$ g'' (x) g^2 (x) = \tfrac34 \beta^3 x^\frac52 $, $x \in \R_+$. 
In a similar way, 
\begin{equation}
    \mathbb{J}_{22}
    \leq
    C
\Big(
1 +  \sup_{s \in [0, T] } \| X ( s ) \|_{ L^{ 5 } ( \Omega; \R^d ) }^{ 5 }
\Big) h^3,
\end{equation}
\begin{equation}
\begin{split}
    \mathbb{J}_{23}
  &  \leq
    6 h
\int_{ t_{j-1} }^{ t_{j} }
\E
[
\|
 g'g (  \widetilde{X} ( t_{j-1} ) ) 
- g'g (  X ( t_{j-1} ) )
\|^2
]
\, \dd  s
\\
 &   \leq
    C
\Big(
1 +  \sup_{s \in [0, T] } \| X ( s ) \|_{ L^{ 6 } ( \Omega; \R^d ) }^{ 6 }
\Big) h^4.
\end{split}
\end{equation}
Altogether we get
\begin{equation}
\mathbb{J}_{2}
    \leq
C h^3 \Big( 1 + \sup_{s \in [0, T]} \| X_s \|_{L^6 ( \Omega; \R) }^6 \Big).
\end{equation}
Inserting this into \eqref{eq:Rj2-decomp-32model} gives
\begin{equation}
\E [ \| \mathcal{R}_{ j } \|^2 ]
\leq
C h^3 \Big( 1 + \sup_{s \in [0, T]} \| X_s \|_{L^6 ( \Omega; \R) }^6 \Big).
\end{equation}
Also,
\begin{equation} \label{Chapter5eq:Ri-conditional-decomposition}
\begin{split}
& \E \big[ \big \| \E ( \mathcal{R}_{ j } | \mathcal{F}_{t_{ j -1} } ) \big \| ^2 \big] 
\\
& \quad = 
\E \Big[ 
\Big\|
\int_{ t_{j-1} }^{ t_{ j } } 
\E \big( 
[ f( X ( s ) ) - f( \widetilde{X} ( t_{j-1} ) ) ]
\vert \mathcal{ F }_{ t_{ j - 1} } 
\big)
\dd s 
\Big\|^2
\Big]
\\
& \quad \leq
2 
\E \Big[ 
\Big\|
\int_{ t_{j-1} }^{ t_{ j } } 
\E \big( 
[ f( X ( s ) ) - f( X ( t_{j-1} ) ) ]
\vert \mathcal{ F }_{ t_{ j - 1} } 
\big)
\dd s 
\Big\|^2
\Big]
\\ & \qquad
+ 
2 
\E \Big[ 
\Big\|
\int_{ t_{j-1} }^{ t_{ j } } 
\E \big( 
[ f( \widetilde{X} ( t_{j-1} ) ) - f( X ( t_{j-1} ) ) ]
\vert \mathcal{ F }_{ t_{ j - 1} } 
\big)
\dd s 
\Big\|^2
\Big].
\end{split}
\end{equation}
Similar to the above arguments, one can show
\begin{equation} \label{Chaper5eq:f-conditional-expect-32-model}
\begin{split}
&
\E \Big[ 
\Big\|
\int_{ t_{j-1} }^{ t_{ j } } 
\E \big( 
[ f( X ( s ) ) - f( X ( t_{j-1} ) ) ]
\vert \mathcal{ F }_{ t_{ j - 1} } 
\big)
\dd s 
\Big\|^2
\Big]
\\ &
\quad =
\E \Big[ 
\Big\|   
\int_{ t_{j-1} }^{ t_{ j } } 
\int_{t_{j-1} }^{s} 
\E
\big[
\big(
f' ( X(r) ) f ( X(r) ) + \tfrac12 f '' ( X(r) ) g^2 ( X(r) )
\big)
\vert \mathcal{ F }_{ t_{ j - 1} } 
\big] 
\, \dd r
\, \dd s 
%\bigg) 
\Big\|^2 \Big] 
\\ &
\quad \leq
C h^4 \Big( 1 + \sup_{s \in [0, T]} \| X (s) \|_{L^6 ( \Omega; \R) }^6 \Big).
\end{split}
\end{equation}
Still, repeating the above arguments shows
\begin{equation}
\begin{split} 
& \E \Big[ 
\Big\|
\int_{ t_{j-1} }^{ t_{ j } } 
\E \big( 
[ f( \widetilde{X} ( t_{j-1} ) ) - f( X ( t_{j-1} ) ) ]
\vert \mathcal{ F }_{ t_{ j - 1} } 
\big)
\dd s 
\Big\|^2
\Big]
\\
& \quad 
\leq
    h^2
\E [
\|
f ( \widetilde{X} ( t_{j-1} ) ) - f (  X ( { t_{j-1} } ) )   \|^2
]
\\
& \quad \leq
C
h^4
\Big(
1 +  \sup_{s \in [0, T] } \| X ( s ) \|_{ L^{ 6 } ( \Omega; \R^d ) }^{ 6 }
\Big).
\end{split}
\end{equation}
These two estimates plugged into \eqref{Chapter5eq:Ri-conditional-decomposition} help us get
\begin{equation}
    \E \big[ \big \| \E ( \mathcal{R}_{ j } | \mathcal{F}_{t_{ j -1} } ) \big \| ^2 \big]
    \leq
    C
h^4
\Big(
1 +  \sup_{s \in [0, T] } \| X ( s ) \|_{ L^{ 6 } ( \Omega; \R^d ) }^{ 6 }
\Big).
\end{equation}
Thanks to Proposition \ref{Chpter5lem:32model-scheme-well-posedness}, the assumption $ \alpha \geq \tfrac52 \beta^2 $ ensures
\begin{equation}\label{eq:32model-L6-MB}
\sup_{ s \in [ 0, T] }
\lVert X( s )  \rVert_{L^6 ( \Omega; \R )} < \infty.
\end{equation}
Now by Theorem \ref{Chaper5Milstein-thm:convergence-rate}
we arrive at the desired assertion
\eqref{eq:Chapter5thm-32model-converg-rate}.
\hfill $\square$

\subsection{Optimal complexity of MLMC}
This subsection is devoted to the analysis of the complexity of MLMC combined with the proposed time-stepping scheme.
\begin{lem}
\label{lem:variance}
Assume the payoff function $\phi : \R \rightarrow \R$ is Lipschitz, i.e.,
\begin{equation}\label{eq:phi-Lipschitz}
| \phi (x) - \phi ( y ) | 
\leq
L_{\phi} | x - y|.
\end{equation}
Let $X_0 > 0$ and $\mu , \alpha , \beta > 0$ 
satisfying $ \alpha \geq \tfrac52 \beta^2 $. 
Let $\{ X({ t}) \}_{ t \in [0, T]}$ be the unique positive solutions to the 3/2 model \eqref{Milstein-eq:SODE}-\eqref{eq:f-g-notation}. 
Let $ Y_{N_{\ell}} $ be produced by 
\eqref{Chpter5Milstein-eq:scheme} using stepsize $h=T/N_{\ell}$.
Let  
\begin{equation}
\bar{\Phi}_{\ell} = \phi ( Y_{ N_{\ell} } )
\end{equation}
and 
%let the final multilevel estimator $\bar{Z}$ be defined by 
%\eqref{equation:final-multilevel-estimator}.
Let $\Phi^X$, $\bar{Z}_{\ell}$ and  $\bar{Z}$ be defined as  \eqref{def:functional-P-in-mlmc}, \eqref{definition:estimator-of-multilevel-estimator} and \eqref{equation:final-multilevel-estimator}, respectively.
Then
\begin{equation}
\text{Var} \left[\bar{Z}_{\ell}\right] \leq K_{2} N_{\ell}^{-1} h_{\ell}^{2}.
\end{equation}
\end{lem}
{\bf Proof}. The assertion for the case $\ell = 0$ is trivial. For $\ell >0$, we combine \eqref{eq:Chapter5thm-32model-converg-rate} with \eqref{eq:phi-Lipschitz} to deduce
\begin{equation}
\begin{split}
\text{Var} \left[\bar{Z}_{\ell}\right]
& =
N_{\ell}^{-2}\sum_{i=1}^{N_{\ell}}
  \text{Var}
  \left[
  \bar{\Phi}_{\ell}^{(i)}-\bar{\Phi}_{\ell -1}^{(i)} 
  \right]
\\
&
=
N_{\ell}^{-1}
  \text{Var}
  \left[
  \bar{\Phi}_{\ell}-\bar{\Phi}_{\ell -1} 
  \right]
\\
&
\leq
N_{\ell}^{-1}
\E 
\left[
  | \phi ( Y_{ N_{\ell} } ) - \phi ( Y_{ N_{\ell-1} } ) |^2 
  \right]
\\
&
\leq
2 N_{\ell}^{-1}
\E \left[
  | \phi ( Y_{ N_{\ell} } )- \phi (X(T) ) |^2 
  \right]
  +
2 N_{\ell}^{-1} \E 
\left[
  |  \phi ( X(T) ) - \phi ( Y_{ N_{\ell-1} } ) |^2 
  \right]
\\
&
\leq
2 N_{\ell}^{-1} 
L_{\phi}^2
h_{ \ell }^2
+
2 N_{\ell}^{-1} 
L_{\phi}^2
h_{ \ell -1 }^2
\\
&
\leq
K_2 N_{\ell}^{-1} 
h_{ \ell }^2,
\end{split}
\end{equation}
as required.
\hfill $\square$
\begin{cor}
\label{cor:optimal-complexity}
Let all conditions in Lemma \ref{lem:variance} be fulfilled.
Then  there exists a positive constant $K_{4}$ such that for any $\epsilon < e^{-1}$ there are values $L$ and $N_{\ell}$ \hspace{0.05em} for which the multilevel estimator \eqref{equation:final-multilevel-estimator}
has a mean-square-error  
\begin{equation}
\text{MSE} := \mathbb{E}\left[(\bar{Z}-\mathbb{E}[\Phi^X])^{2}\right]
<\epsilon^{2}  \notag
\end{equation}
with a  computational complexity $\mathcal{C} \leq K_4 \epsilon^{ - 2 }$.

{\bf Proof}. 
Thanks to \eqref{eq:Chapter5thm-32model-converg-rate}
and \eqref{eq:phi-Lipschitz} and by the H\"{o}lder inequality, 
the condition (a) in Theorem \ref{theorem:mlmc-complexity-theorem} is fulfilled with $\chi = 1$. 
Conditions (b) and (d) are trivially satisfied,
with $\theta = 1$.
Further, the condition (c)
has been validated by Lemma \ref{lem:variance} with $\gamma = 2 > \theta$. 
As a consequence of \eqref{eq:comlexity-thm-C}, the computational complexity $\mathcal{C} \leq K_4 \epsilon^{ - 2 }$ follows.
\hfill $\square$

\end{cor}

\section{Numerical experiments}
\label{Milstein-sect:numer-results}

\begin{figure}
\centering
      \includegraphics[width=5in,height=4in]  {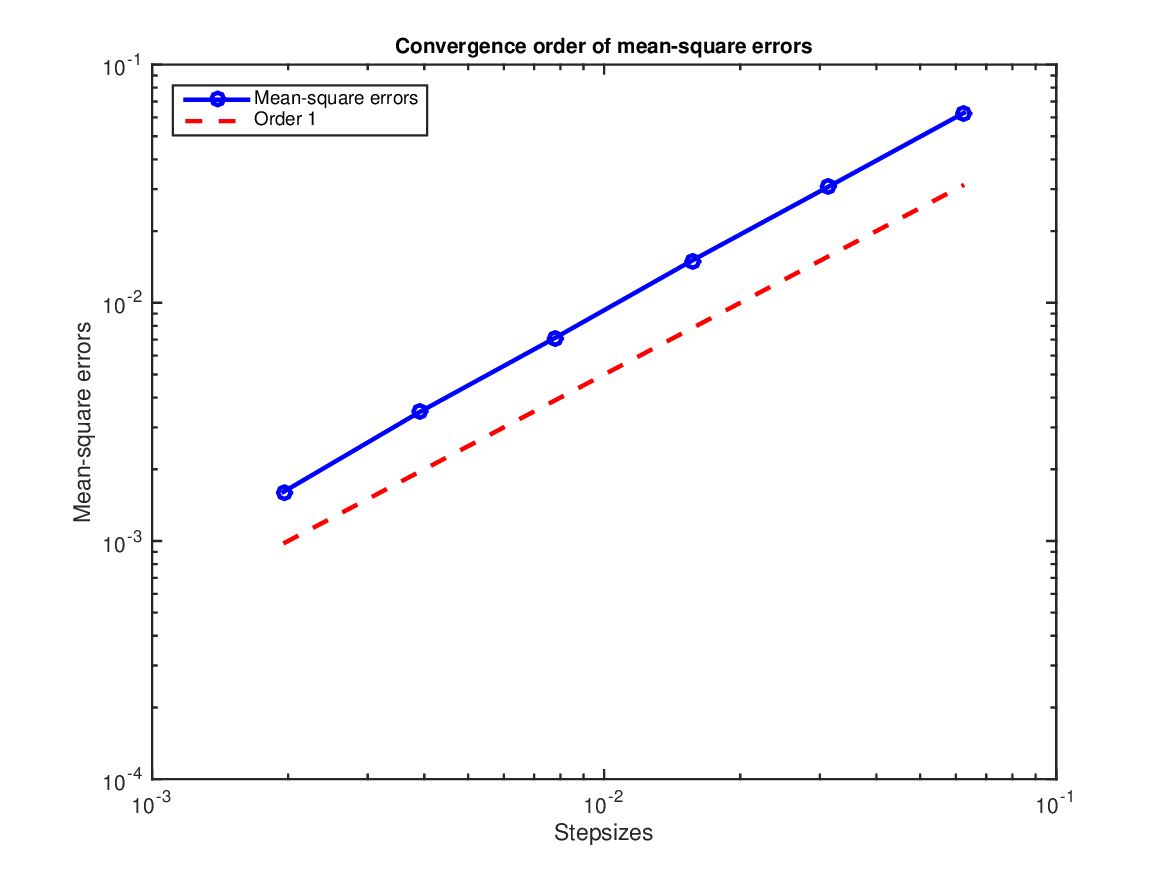}
 \caption{Mean-square convergence rates of the proposed time-stepping scheme.}
\label{Milstein-fig:convergence-rates-32model}
\end{figure}

This section provides some numerical simulations to support the above theoretical findings.
Let us consider the 3/2-model
\begin{equation}
 d X (t) = X(t) ( \mu - \alpha X(t) ) dt + \beta X(t) ^ {3/2} \, \dd W(t), \quad
 \mu , \alpha , \beta > 0,
 \quad
     X(0) = 1,
\end{equation}
%. 
where  we take $\mu = 1, \alpha = 5/2, \beta = 1$ for the following numerical tests. Evidently, $ \alpha \geq \tfrac52 \beta^2 $ in our setting and all the above conditions are fulfilled.
As usual, the expectations are approximated by computing averages over $5000$ samples and the ``exact'' solutions are identified as numerical approximations using a fine stepsize $h_{\text{exact}} = 2^{-12}$.
%To confirm the convergence rates of SST methods, 
In Figure \ref{Milstein-fig:convergence-rates-32model}  we depict mean-square errors of the proposed scheme using six different stepsizes $h = 2^{-i}, i =4,5,...9$.
From Figure \ref{Milstein-fig:convergence-rates-32model},
one can easily see that the resulting approximation errors 
shrink at a slope very close to $1$, 
which well agrees with the previously obtained convergence rate.
%\noindent $\bullet$ Case I:  $\kappa = 4, \, \rho = 2, ...$
%\vspace{0.5cm}
%
%\noindent $\bullet$ Case II:   $\kappa = 3, \, \rho = 2, ...$
%\newline

%

\begin{figure}
\centering
      \includegraphics[width=5in,height=4in]{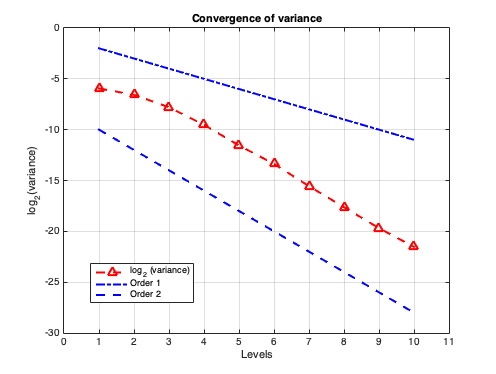}
 \caption{Convergence of variance}
\label{fig:variance-picture}
\end{figure}

\begin{figure}
\centering
      \includegraphics[width=5in,height=4in]{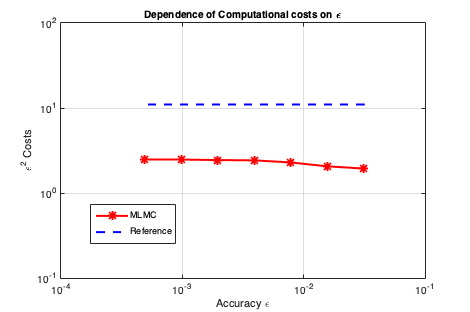}
 \caption{Computational costs}
\label{fig:complexity-picture}
\end{figure}

Next we test the convergence of the variance as a function of the level of approximation and the computational cost for the payoff function $\phi (x) = \max(x - 0.05, 0)$. We run $10^{5}$ samples to compute the variance of the multilevel estimators at $10$ different levels in Figure \ref{fig:variance-picture},
where reference lines with slopes of order $-1$, $-2$ are also added. It can be easily observed that 
$\gamma = 2$ in the MLMC setting.
In  Figure \ref{fig:complexity-picture}, the dependence of complexity as function of the desired accuracy $\epsilon$ for the MLMC is presented. 
% for $ h = 2^{-5}, 2^{-6}..., 2^{-11}$
There one can clearly see 
that the complexity of the MLMC is proportional to $\epsilon^{2}$, which confirms the above theoretical findings.

%\begin{table}[!ht]
%\label{Milstein-table:least-squares}
%\begin{center}% \footnotesize
%\caption{A least squares fit for the convergence rate $\delta$.} 
%\begin{tabular*}{12cm}{@{\extracolsep{\fill}}ccc}
%\hline  & Case I  &  Case II
% \\ \hline
%$\theta=\frac12$ & $\delta$ = 0.5588, resid = 0.0554  & $\delta$ = 0.5571, resid = 0.0499
%\\ \hline
%$\theta=1$ & $\delta$ = 0.5554, resid = 0.0546 & $\delta$ = 0.5633, resid = 0.0396   \\
% \hline
%\end{tabular*}
%\end{center}
%\end{table}

\section{Conclusion}
\label{sec:conclusion}

In this work we are interested in quantification of expectations of some functions of the solution $X_{T}$ to
the Heston 3/2-model. To this end, we turn to  
the MLMC approach combined with a newly proposed 
Milstein-type scheme for the Heston 3/2-model. 
It is shown that the time-stepping scheme is
explicit and positivity preserving for any stepsize $h>0$,
which is particularly desirable for large discretization time steps naturally appearing in the MLMC setting.  
Moreover, we prove the mean-square convergence rate
of order one for the scheme applied to the 3/2-model.
This in turn promises the desired relevant variance of the multilevel estimator
        and thus justifies the optimal complexity $\mathcal{O}(\epsilon^{-2})$ for the MLMC approach.
As a possible future work, we aim to do the analysis
for more SDE models in practice \cite{neuenkirch2014first} and for 
options with non-globally Lipschitz payoff
\cite{giles2009analysing}.

%\bibliographystyle{abbrv}

%\bibliography{../bib/bibfile}
%\bibliography{bibfile}

\end{document}